\documentclass[11pt]{article}

\usepackage{setspace}
\usepackage{graphicx}
\usepackage{amsmath}
\usepackage{amsfonts}
\usepackage{epsfig}
\usepackage{epstopdf}
\usepackage{color}
\newcommand{\be}{\begin{equation}}
\newcommand{\ee}{\end{equation}}
\newcommand{\bes}{\begin{equation*}}
\newcommand{\ees}{\end{equation*}}
\newcommand{\beqn}{\begin{eqnarray}}
\newcommand{\eeqn}{\end{eqnarray}}
\newcommand{\beqns}{\begin{eqnarray*}}
\newcommand{\eeqns}{\end{eqnarray*}}

\newcommand{\lkr}{\left(}
\newcommand{\lkv}{\left[}
\newcommand{\rkv}{\right]}
\newcommand{\rkr}{\right)}
\newcommand{\lfi}{\left\{}
\newcommand{\rfi}{\right\}}

\newcommand{\fr}[1]{(\ref{#1})}

\newcommand{\ph}{\varphi}
\newcommand{\del}{\delta}

\newcommand{\af}{\alpha}
\newcommand{\eps}{\epsilon}

\newcommand{\te}{\theta}
\newcommand{\om}{\omega}
\newcommand{\lam}{\lambda}
\newcommand{\Up}{\Upsilon}

\newcommand{\sig}{\sigma}

\newcommand{\Om}{\Omega}

\newcommand{\PP}{\ensuremath{{\mathbb P}}}

\newcommand{\RR}{\ensuremath{{\mathbb R}}}

\newcommand{\Span}{\mbox{Span}}

\newcommand{\Cov}{\mbox{Cov}}
\newcommand{\diag}{\mbox{diag}}

\newcommand{\Tr}{\mbox{Tr}}
\newcommand{\proj}{\mbox{proj}}

\newtheorem{theorem}{Theorem}
\newtheorem{lemma}{Lemma}
\newtheorem{corollary}{Corollary}

\newcommand{\bd}{\mathbf{d}}

\newcommand{\bof}{\mathbf{f}}

\newcommand{\bq}{\mathbf{q}}
\newcommand{\bt}{\mathbf{t}}
\newcommand{\bu}{\mathbf{u}}
\newcommand{\bv}{\mathbf{v}}

\newcommand{\bx}{\mathbf{x}}
\newcommand{\by}{\mathbf{y}}
\newcommand{\bz}{\mathbf{z}}

\newcommand{\bA}{\mathbf{A}}

\newcommand{\bD}{\mathbf{D}}

\newcommand{\bQ}{\mathbf{Q}}

\newcommand{\bU}{\mathbf{U}}

\newcommand{\bW}{\mathbf{W}}
\newcommand{\bX}{\mathbf{X}}

\newcommand{\bte}{\mbox{\mathversion{bold}$\te$}}

 \newcommand{\btheta}{\mbox{\mathversion{bold}$\theta$}}

\newcommand{\boeta}{\mbox{\mathversion{bold}$\eta$}}

\newcommand{\bzeta}{\mbox{\mathversion{bold}$\zeta$}}
\newcommand{\bphi}{\mbox{\mathversion{bold}$\ph$}}
\newcommand{\bpsi}{\mbox{\mathversion{bold}$\psi$}}

 \newcommand{\hbtheta}{\widehat{\btheta}}

\newcommand{\bPhi}{\mbox{\mathversion{bold}$\Phi$}}
\newcommand{\bUp}{\mbox{\mathversion{bold}$\Up$}}
\newcommand{\bPsi}{\mbox{\mathversion{bold}$\Psi$}}
\newcommand{\bLam}{\mbox{\mathversion{bold}$\Lambda$}}
\newcommand{\bOm}{\mbox{\mathversion{bold}$\Om$}}

\newcommand{\Jc}{{J_{*}^c}}

\newcommand{\calJ}{{\mathcal{J}}}

\newcommand{\calL}{{\mathcal{L}}}
\newcommand{\calP}{{\mathcal{P}}}

\newcommand{\calN}{{\cal N}}

\newcommand{\hbte}{\widehat{\bte}}
\newcommand{\hte}{\widehat{\te}}

\newcommand{\hbf}{\hat{\bof}}
\newcommand{\hbq}{\hat{\bq}}

\newcommand{\hk}{\widehat{k}}

\newcommand{\alfo}{\af_0}

\newcommand{\sumjp}{\sum_{j=1}^p}

\newcommand{\sumJ}{\sum_{j \in J_{*}}}

\newcommand{\bphij}{\bphi_j}

\newcommand{\nuj}{\nu_j}

\newcommand{\lamin}{\lam_{\min}}
\newcommand{\lamax}{\lam_{\max}}

\begin{document}

\title{\bf { Solution of linear   ill-posed problems using  random dictionaries }}

\author{{\em  Pawan Gupta and  Marianna Pensky}   \\
         Department of Mathematics,
         University of Central Florida      }

\date{}

\bibliographystyle{plain}
\maketitle

\begin{abstract}
In the present paper we consider application of  overcomplete dictionaries to solution of general ill-posed linear inverse problems. 
In the context of regression problems, there has been enormous amount of effort to recover an unknown function
using such dictionaries. One of the most popular methods, lasso and its versions, is based 
on minimizing empirical likelihood and unfortunately, requires 
stringent assumptions on the dictionary, the, so called, compatibility conditions.  
Though compatibility conditions are hard to satisfy, it is well known that this can  be accomplished by using random dictionaries.
In the present paper, we show how one can apply random dictionaries to solution of  ill-posed 
linear inverse problems. We put a theoretical foundation under the suggested methodology and study its performance
via simulations.

 \vspace{2mm} 


 \vspace{2mm}
{\bf AMS (2000) Subject Classification}: {Primary: 62G05. Secondary: 62C10    }
\end{abstract}

\section{Introduction  }
\label{sec:introduction}
\setcounter{equation}{0}

In this paper, we consider   solution of a general ill-posed linear inverse problem $\bQ \bof =  \bq$ where $\bQ$ is a  
bounded linear  operator that  does not have a bounded inverse and the right-hand side $\bq$ is measured with error. 
In particular, we consider equation 
\be \label{geneq}
\by  =  \bq  +  \sigma \boeta, \quad \bq  =  \bQ \bof,  
\ee
where $\by,  \bq, \bof, \boeta \in \RR^n$, $\bQ \in \RR^{n \times n}$. Here, $\by$ is observed,  $\bq$ is unobserved, $\bof$ 
is the function to be estimated,     $\sig$ is the noise level and $\boeta$ is the noise vector 
which we assume to have  the standard normal distribution.  Matrix $\bQ$ is invertible but its lowest eigenvalue is very small, 
especially, when $n$ is relatively large,   which makes the problem ill-posed.   A general linear inverse problem can usually be reduced to 
formulation  \fr{geneq} by  either  expanding  $\by$ and $\bof$ 
over some collection of basis functions  or by measuring them at some set of points.

Solutions of statistical inverse problem  \fr{geneq}   usually rely on reduction of the problem 
to the sequence model by   carrying out the singular value decomposition (SVD)  (see, e.g.,   \cite{cavalgol1},
\cite{cavreiss},   and \cite{tropp} and references therein), or its relaxed version, the 
wavelet-vaguelette decomposition proposed by Donoho \cite{donoho} and further studies by Abramovich and Silverman \cite{abr}. 
Another general approach  is Galerkin method with subsequent model selection  (see, e.g.,    
 \cite{cohen}).

The advantage of the methodologies listed above is that they are asymptotically optimal in a minimax 
sense.  The function of interest is usually represented via an   orthonormal basis which is
motivated by the form of matrix $\bQ$. However, in spite of being minimax optimal in many contexts, 
these approaches   have  drawbacks. In particular, in practical applications, the number of observations $n$ may be low 
while   noise level $\sigma$  high. In this situation, if the unknown vector $\bof$ does not 
have a relatively compact and accurate representation in the chosen basis,  
the precision of the resulting estimator will be  poor.

In the last decade, a great deal of effort was spent on recovery of an unknown vector $\bof$ in 
regression setting from its noisy observations using overcomplete  dictionaries. 
In particular, if   $\bof$ has a sparse representation in some dictionary  
(a collection of vectors used for the representation of $\bof$), 
then $\bof$ can be recovered with a much better precision than, for example,  when  it is expanded over an orthonormal basis. 
The methodology is based on the idea  that the error of an estimator of $\bof$  is approximately
proportional  to the number of dictionary functions that are used for   representing   $\bof$,  
therefore,  expanding a function of interest over fewer dictionary elements reduces 
the estimation error.  Similar advantages hold in the case of linear inverse problems (see \cite{pensky}).
However, in order to represent a variety of functions using a small number 
of dictionary elements, one   needs to consider a dictionary of much larger 
size than the number of available  observations, the, so called, overcomplete dictionary.

A variety of techniques have been developed for solution of   regression problems using overcomplete dictionaries
% those problems 
such as likelihood  
penalization methods and greedy algorithms. The most popular of those methods (due to its computational convenience),  
lasso and its versions,  have been used  for solution of a number of theoretical and applied statistical 
problems (see, e.g.,  \cite{bickel_ritov_tsybakov},   %\cite{bunea_tsybakov_2},  \cite{arnak},
% \cite{lounici_pontil_tsybakov},   
and also \cite{sara} and references therein).
However, application of lasso is based on maximizing the   likelihood and, unfortunately, 
 relies on stringent assumptions on the dictionary  $\lfi \bphi_k \rfi_{k=1}^p$, the, so called, compatibility conditions,
for a proof of its optimality.  
In regression set up, as long as compatibility conditions hold, lasso identifies a linear combination of   
the dictionary elements which represent the function of interest 
best of all at a ''price'' which is proportional to $\sigma  \sqrt{n^{-1} \, \log p}$  where $p$ is the dictionary size
(see, e.g., \cite{sara}). Regrettably, while compatibility conditions may be satisfied for the vectors $\bphi_j$ in the  original dictionary,
they usually do not hold for their images $\bQ \bphi_j$  due to contraction imposed by the   operator $\bQ$.
Pensky \cite{pensky} showed how lasso solution can be modified, so that it delivers an  optimal solution,
however, compatibility assumptions in \cite{pensky} remain very complex and hard to verify.

In the recent years it has been discovered  that in regression setting, one can  satisfy compatibility conditions  for lasso
by simply using random dictionaries. In particular, Vershynin \cite{vershynin}  provided a variety of way 
for construction of such dictionaries, i.e, dictionaries comprised of random vectors.  
The methodology  of \cite{vershynin}, however, is intended for the recovery of a function which is directly observed. 
The purpose of the present paper is to explain how random dictionaries can be adopted for solution of 
ill-posed linear inverse problems.

To the best of our knowledge, application of random dictionaries to signal recovery has not been attempted so far since 
random vectors usually contain ``pure noise'' and therefore are perceived as unsuitable for representing 
a meaningful signal. This is indeed true when one needs to estimate a  simple smooth function which is best represented by a small 
set of smooth basis functions. However, when a signal has a more complicated structure, it cannot be expanded over a small number 
of basis functions. In this case, a  large rich  dictionary may be helpful since there is a high chance that the signal of interest
can be well approximated by a linear combination of a small number of  vectors of the dictionary. The advantage of the random dictionaries is that, 
unlike in the case of fixed dictionaries, one can work with a dictionary of vary large size which provides a competitive advantage 
over standard orthogonal basis based techniques. This benefits can be more significant when one needs to solve 
an inverse ill-posed problem since, as it was notes in \cite{pensky},
 finding a ``low-cost'' representation of a function of interest can significantly improve the accuracy of the solution.

The rest of the paper is organized as follows. Section~\ref{sec:method} introduces some notations,
formulates optimization problem with lasso penalty and lists compatibility conditions of \cite{pensky}.
Section~\ref{sec:random} contains the main results of the paper: it explains how one can obtain fast lasso convergence rates 
by using random dictionaries. Section~\ref{sec:examples} contains a simulation study which proves that our
technique is competitive. Section~\ref{sec:discussion} concludes the paper with the discussion.
Finally, Section \ref{sec:proofs} contain proofs of the statements in the paper.

%%%%%%%%%%%%%%%%%%%%%%%%%%%%%%%%%%%%%%%%%%%%%%%%%%%%%%%%%%%%%%%%%%%%%%%%%%%%%%%%%%%%%%%%%%%%%%%%%%%%%%%%%%%%%%%%%%%%55

\section{Construction of the lasso estimator and a general compatibility condition}
\label{sec:method}
\setcounter{equation}{0}

In the paper, we use the following notations.

For any vector $\bt \in \RR^p$, denote  its $\ell_2$, $\ell_1$, $\ell_0$ and $\ell_\infty$ norms by, 
respectively,  $\| \bt\|$, $\| \bt\|_1$,  $\| \bt\|_0$ and $\| \bt\|_\infty$.
For any matix $\bA$, denote its $i^{th}$ row and $j^{th}$ column by, $\bA_{i \cdot}$ and $\bA_{ \cdot j}$ respectively.
Denote its spectral and Frobenius norms by, respectively,  $\| \bA \|$ and $\| \bA \|_2$.
% Notation $\bA >0$ or $\bA \geq 0$ means, respectively,  that $\bA$ is positive or non-negative definite.
% Denote determinant of $\bA$ by $|\bA|$ and the largest, in absolute value, element of $\bA$ by $\| \bA\|_{\infty}$.
% Denote the Moore-Penrose inverse of matrix $\bA$ by $\bA^{+}$.
% 
Denote $\calP = \{1, \cdots, p\}$.   For any subset of indices  $J \subseteq \calP$, 
 subset $J^c$ is its complement in $\calP$ and  $|J|$ is its cardinality, so that $|\calP| =p$. 
Let  $ \calL_J = \Span \lfi \bphi_j, \   j \in J \rfi$. 
% 
% \item
If $J \subset \calP$ and $\bt \in \RR^p$,  then     $\bt_J \in \RR^{|J|}$ denotes   reduction of vector $\bt$ to
subset of indices $J$. Also, $\bPhi_J$ denotes the reduction of matrix $\bPhi$ to columns $\bPhi_{ \cdot j}$
with $j \in J$.

Denote by $\lamin (m; \bPhi)$  and $\lamax (m; \bPhi)$ the minimum and the maximum restricted eigenvalues
of matrix $\bPhi^T \bPhi$  given by  
\be \label{eigrestrict}
\lamin (m; \bPhi) = \min_{\stackrel{\bt \in \RR^p}{\|\bt \|_0 \leq m}}\ \frac{\bt^T \bPhi^T \bPhi \bt}{\| \bt \|_2^2}, 
\quad
\lamax (m; \bPhi) = \max_{\stackrel{\bt \in \RR^p}{\|\bt \|_0 \leq m}}\ \frac{\bt^T \bPhi^T \bPhi \bt}{\| \bt \|_2^2}. 
\ee
% Whenever there is no ambiguity, we drop $\bPhi$ in the above notations and write simply  
% $\lamin (m)$ and $\lamax (m)$.

% \item
%  $a_m \asymp b_m$ means that there exist    constants $0 < C_1<C_2<\infty$
% independent of $m$ such that $C_1 a_m < b_m < C_2 a_m$.
%\end{itemize}
%

Denote by $\bPhi$ the dictionary matrix with columns $\bphij \in \RR^n$, $j=1, \cdots, p$, where $p$ is possibly much larger than  $n$ 
and 
 \be \label{f_expan}
\bof_{\bt} = \sum_{j=1}^p  t_j  \bphi_j = \bPhi \bt.
\ee
Let $\btheta$ be the true vector of coefficients of expansion of $\bof$ over the dictionary $\bPhi$,
so that $\bof = \bPhi \btheta$. Let vectors $\bpsi_j$ be such that $\bQ^T \bpsi_j = \bphi_j$, where 
$\bQ^T$ is the   transpose of matrix $\bQ$, and $\bPsi$ be a matrix  with columns   $\bpsi_j$, $j=1, \cdots, p$.
Then, 
\be   \label{Psi}
 \bQ^T \bPsi = \bPhi \quad \mbox{and} \quad \bPsi = \bQ (\bQ^T \bQ)^{-1}   \bPhi. 
\ee 
Note that, although $\bof$ is unknown,  
  \be  \label{ernorm}
 \| \bof - \bof_{\bt} \|^2 = \| \bof \|_2^2 + \bt^T \bPhi^T \bPhi \bt  - 2 \bt^T \bPhi^T \bof    =
\| \bof \|_2^2 + \bt^T \bPhi^T \bPhi \bt  - 2 \bt^T \bPsi^T \bQ \bof   
 \ee 
is the sum of the three components where the first one  is independent of $\bt$, 
 the second one  is completely known, while the last term is of the form $2 \bt^T \bPsi^T\bQ \bof = 2 \bt^T \bPsi^T \bq $ and, hence,
can   be estimated by $2 \bt^T \bPsi^T \by$.   Let $\bz$ be such that
\bes
\bPsi^T \by = \bPhi^T \bz. 
\ees
Therefore,  expression $\bt^T \bPhi^T \bPhi \bt  - 2 \bt^T \bPsi^T \by$ is minimized by the same vector $\bt$ that minimizes
$\| \bPhi \bt - \bz \|_2^2$ where
\be \label{bz}
\bz = (\bPhi \bPhi^T)^{-1} \bPhi \bPsi^T \by.
\ee

Denote $\nu_j = \|\bpsi_j\|_2$, $j=1, \cdots, p$, and  observe that $\nu_j$ is proportional to  the standard deviation of the 
$j$-th component of the vector  $\bPsi^T  \by$. The value  $\nu_j$ can be viewed as a ``cost'' of using a dictionary element $\bphi_j$
in representation of $\bof$. Consider a  matrix 
\be \label{gamma}
\bUp = \diag(\nu_1, \cdots, \nu_p)=\diag(\|\bpsi_1\|_2, \cdots, \|\bpsi_p\|_2). 
\ee
Following \cite{pensky}, we estimate  the true vector of coefficients $\btheta$ as a solution of 
the quadratic optimization problem with the weighted lasso penalty
\be \label{las_sol}
\widehat{\btheta} = \arg\min_{\bt}    \lfi  \| \bPhi \bt - \bz \|_2^2  + 
\af    \| \bUp \bt\|_1   \rfi.
\ee
Here $\bz$ is given by \fr{bz} and  $\af \geq \alfo$ where
\be \label{alf0}
\alfo =  \sigma\ \sqrt{2\,  n^{-1} \, (\tau +1) \log p}.
\ee
Parameter $\tau >0$ is  related to the required probability bound (see formula \fr{eq:setOm0}  in Section~\ref{sec:proofs}
for details).
Subsequently, we estimate the unknown solution $\bof$ by $\widehat{\bof} = \bPhi \widehat{\btheta}$.

Note   that since we are interested in $\bof_{\btheta}$ rather than in the vector $\btheta$ of coefficients 
themselves, we are using lasso  for solution of the  so called  prediction problem where it requires 
milder conditions on the dictionary. In fact, it is known (see \cite{pensky})
that with no additional assumptions,  for  $\af \geq \alfo$, with probability at least
$1 - 2 p^{-\tau}$,  one has
\be \label{slowlas}
n^{-1}\, \| \bof_{\hbtheta} - \bof \|_2^2 \leq \inf_{\bt } \lkv  n^{-1}\, \| \bof_{\bt } - \bof \|_2^2 + 4 \af   
\sumjp \nuj |t_j|  \rkv.
\ee
It is easy to see that if $\bt = \btheta$, then $\bof_{\bt } = \bof$. Then,  with high probability, the error of the 
estimator $\bof_{\hbtheta}$ is proportional to 
$ 
  \sigma \sqrt{n^{-1} \, (\tau +1) \log p}\ \sum_j \nu_j.
$ 
This is the, so called, {\it slow lasso rate}. In order to attain the {\it fast lasso rate}   
proportional to $\sigma^2 n^{-1} \sum_j \nu_j^2$, one needs some kind of a  compatibility assumption.

Pensky \cite{pensky} formulated the following compatibility condition:  matrices $\bPhi$ and $\bUp$ 
are such that for some $\mu >1$ and any $J \subset \calP$
\be \label{comp}
\kappa^2 (\mu, J) = \min \lfi \bd \in \calJ(\mu, J),\, \| \bd \|_2 \neq 0: \quad 
\frac{\bd^T \bPhi^T \bPhi \bd \cdot \Tr(\bUp_J^2)}{\|(\bUp \bd)_{J}\|_1^2} \rfi  >0.
\ee 
where $\calJ (\mu, J)  = \lfi \bd \in \RR^p:\ \|(\bUp \bd)_{\Jc}\|_1 \leq \mu  \|(\bUp \bd)_{J}\|_1  \rfi$.
Pensky \cite{pensky} proved that, under assumption \fr{comp}, for $\af = \varpi \alfo$ where 
$ \varpi \geq  (\mu +1)/(\mu -1)$ and $\af_0$ is  defined in \fr{alf0}, with probability at least $1 - 2 p^{-\tau}$,  one has
\be \label{fasrlas}
\| f_{\hbte} - f \|_2^2   \leq \inf_{J \subseteq \calP} \lkv \| f - f_{\calL _J} \|_2^2   
+   \frac{\sig^2  K_0 (1 + \varpi)^2 (\tau +1)}{\kappa^2 (\mu, J)} \, \frac{\log p}{n} \, \sum_{j \in J} \nu_j^2 \rkv,
\ee 
where $f_{\calL _J}  = \proj_{\calL_J} f$.

Note, however,  that  unless matrix $\bPhi$ has orthonormal columns,  assumption \fr{comp} is hard not only 
to satisfy but even to verify since it requires checking it for every subset $J$ in $\calP$. Indeed, sufficient conditions 
listed in Appendix A1 of \cite{pensky} rely on the results of  Bickel {\it et al.} \cite{bickel_ritov_tsybakov} and require 
very stringent conditions on $\lamin (m; \bPhi)$ and entries  $\bUp$ in \fr{gamma}.
In the present paper, we offer an alternative to this approach.

%%%%%%%%%%%%%%%%%%%%%%%%%%%%%%%%%%%%%%%%%%%%%%%%%%%%%%%%%%%%%%%%%%%%%%%%%%%%%%%%%%%%%%%%%%%%%%%%%%%%%%%%%%%%%%%%%%%%55

\section{Lasso solution to linear inverse problems using   random dictionaries} 
\label{sec:random}
\setcounter{equation}{0}

An advantage of using random dictionary lies in the fact that one can ensure, with a high probability,
that the dictionary satisfies a restricted isometry condition (see, e.g.,  \cite{candez} or \cite{foucart}).
In particular, if matrix $\bPhi \in \RR^{n\times p}$ satisfies the restricted isometry property of order $s \geq 1$, then 
$\lamin (s; \bPhi) >0$.  
The latter allows one to formulate the following results.

\begin{theorem}  \label{th:randLasso}
Let  $\bte$ be the  solution of optimization problem \fr{las_sol} with $\af \ge   \alfo$ 
where $\af_0$ is  defined in \fr{alf0}.
Let $\bPhi \in \RR^{n \times p}$ be a random dictionary independent of $\by$ in \fr{geneq}.
% and $\del \in (0,1)$ and  $\rho >0$ be fixed constants.
Denote % $\hat{s} = \|\hbtheta\|_0$ and 
\be \label{Jstar}
 J_{*} = \arg\min \lfi J \subset \calP:\  n^{-1}\,  \| \bof - \bof_{\calL _J} \|_2^2 +      
K_0 \af^2 \ \sum_{j \in J} \nu_j^2 \rfi,
\ee
where $\bof_{\calL_J}  = \proj_{\calL_J} \bof$ and assume that 
$\bPhi$ is such that for some $s$, $1 \leq s \leq n/2$ and  $\del, \eps_1, \eps_2, \eps_3 \in (0,1)$,  the following conditions hold 
\beqn 
\PP \lkr \lamin(2s; \bPhi)   \geq    1-\delta \rkr  & \geq & 1- \eps_1, % \quad    1 \leq s \leq n/2, 
\label{con1}\\
 \PP \lkr | J_{*}| \leq s \rkr & \geq & 1 - \eps_2,
\label{con2}\\
 \PP \lkr \|\hbtheta\|_0 \leq s \rkr & \geq & 1 - \eps_3,
\label{con3}
\eeqn
If $K_0 \geq 4/(1 -\del)^2$ in \fr{Jstar}, then  
\be \label{newfast1}
\PP \lkr \frac{1}{n}\,\| \bof_{\hbtheta} - \bof \|_2^2   \leq \inf_{J \subseteq \calP} \ 
\lkv \frac{1}{n}\,\| \bof - \bof_{\calL _J} \|_2^2 + K_0 \af^2 \ 
% K_0  \,  \frac{\sigma^2  \log p}{n} \ 
\sum_{j \in J} \nu_j^2 \rkv \rkr \geq 1 - 2 p^{-\tau} - \eps_1 - \eps_2 - \eps_3.
\ee 
\end{theorem}
\medskip

Note that for $\af = \af_0$ and $K_0 = 4/(1 -\del)^2$, under assumptions \fr{con1} -- \fr{con3}, formula \fr{newfast1}
yields the following result
\be \label{newfast2}
\PP \lkr \frac{1}{n}\,\| \bof_{\hbtheta} - \bof \|_2^2   \leq \inf_{J \subseteq \calP} \ 
\lfi \frac{1}{n}\,\| \bof - \bof_{\calL _J} \|_2^2 +   
\frac{4\, \sigma^2}{n\, (1 -\del)^2} \   
\sum_{j \in J} \nu_j^2 \rfi \rkr \geq 1 - 2 p^{-\tau} - \eps_1 - \eps_2 - \eps_3.
\ee

%%%%%%%%%%%%%%%%%%%%%%%%%%%%%%%%%%%%%%%%%%%%%%%%%%%%%%%%%%%%%%%%%%%%%%%%%%%%%%%%%%%%%%%%%%%%%%%%%%%%%%%%%%%%%%%%%%%%55
\noindent
As Lemma \ref{lem1} below shows,   assumption \fr{con1} can be guaranteed by choosing a dictionary of a 
particular type.

\begin{lemma}  \label{lem1}
Let matrix $\bPhi \in \RR^{n \times p}$ be independent of $\by$ and satisfy one of the following conditions:\\
a) Matrix $\bPhi$ has independent sub-gaussian isotropic random rows;\\
b) Matrix $\bPhi$ has independent sub-gaussian isotropic random
columns with unit norms;\\
c) Matrix $\bPhi$  is obtained as  
$\bPhi  =(c\,\sqrt{n})^{-1} \bD \bW$
where $\bW\in \RR^{m \times p}$ is a matrix with i.i.d. standard Gaussian entries 
and columns of the matrix $\bD \in \RR^{n \times m}$ form a non-random $c$-tight frame, 
so that for any vector $\bx$, one has $\bx^T\bD\bD^T\bx=c^2 \|\bx\|^2$.

If, for some $\delta\in (0,1)$ and $1 \leq s \leq n/2$,  one has 
 \be  \label{ncond}
 n \geq C_1 \delta^{-2} s \log(e p/s),  
 \ee 
then condition \fr{con1} holds with $\eps_1 \leq   2 \exp(- C_2 \del^2 n)$.
Here, $C_1$ and $C_2$ depend  on the kind of sub-gaussian variables that are involved in formation of $\bPhi$
and are independent of $n$, $m$, $p$, $s$ and $\delta$.   
\end{lemma}

Finally, conditions \fr{con2} and \fr{con3} can be ensured by restricting the set of solutions $\bt$ to 
vectors with cardinality at most $s$. In this case, $\eps_2 = \eps_3 =0$ and the following corollary 
of Theorem \ref{th:randLasso} is valid.

\begin{corollary}  \label{cor:randLasso}
Let  $\bte$ be the  solution of optimization problem  
\be \label{restrict_sol}
\widehat{\btheta} = \arg\min_{\bt: \|\bt\|_0 \leq s}    \lfi  \| \bPhi \bt - \bz \|_2^2  + 
\af    \| \bUp \bt\|_1   \rfi,
\ee
with $\af \ge   \alfo$ where $\af_0$ is  defined in \fr{alf0}.
Let $\bPhi \in \RR^{n \times p}$ be one of the random  random dictionaries defined in Lemma \ref{lem1}.
If, for some $\delta\in (0,1)$,  condition \fr{ncond} holds,
then 
\be \label{newfast3}
\PP \lkr \frac{1}{n}\,\| \bof_{\hbtheta} - \bof \|_2^2   \leq \inf_{\stackrel{J \subseteq \calP}{|J| \leq s}}  
\lkv \frac{1}{n}\,\| \bof - \bof_{\calL _J} \|_2^2 + \frac{4   \af^2}{(1-\del)^2}  
\sum_{j \in J} \nu_j^2 \rkv \rkr \geq 1 - 2 p^{-\tau} - 2 \exp(- C_2 \del^2 n),
\ee 
 where $C_2$ depends  on the kind of sub-gaussian variables that are involved in formation of $\bPhi$
and is independent of $n$, $m$, $p$, $s$ and $\delta$.
\end{corollary}

%%%%%%%%%%%%%%%%%%%%%%%%%%%%%%%%%%%%%%%%%%%%%%%%%%%%%%%%%%%%%%%%%%%%%%%%%%%%%%%%%%%%%%%%%%%%%%%%%%%%%%%%%%%%%%%%%%%%55

Note that case c) above offers a structured random dictionary since each of its elements 
is a linear combination of smooth functions.

%%%%%%%%%%%%%%%%%%%%%%%%%%%%%%%%%%%%%%%%%%%%%%%%%%%%%%%%%%%%%%%%%%%%%%%%%%%%%%%%%%%%%%%%%%%%%%%%%%%%%%%%%%%%%%%%%%%%%%%%%%%%%%%%%%%%%%%%%%%%%%%%%%%

\section{Simulation studies}
\label{sec:examples}
\setcounter{equation}{0}

In order to evaluate the performance of the procedure suggested in this paper we carried out a limited 
simulation study. For our study, we chose three sample sizes $n=32$, $n=64$  and $n=128$.  
% We also chose the number of dictionary elements in the random dictionary as $p=5000$. 
We first generated a true vector $\bof$ using {\tt MakeSignal} 
program in the package Wavelab~850. We then generated the matrix $\bQ$ in  \eqref{geneq} as $\bQ = \bU \bLam \bU^T$ where $\bU$ is an  
 $(n \times n)$ random orthogonal matrix and $\bLam$ is a diagonal matrix with entries $\bLam_{ii} = 1/\sqrt{i}$, $i=1, 2,  \cdots, n$.
Using $\bQ$ we obtained the unobserved vector  $\bq$ as
\bes
\bq=\bQ \bof.
\ees
At last, for generating the   data $\by$ we added Gaussian random noise  to $\bq$.
% Data are always observed, otherwise they are not data
For this purpose,  we chose particular values of the Signal to Noise Ratio (SNR) and obtained $\sigma$ 
as the ratio of the standard deviation of $\bq$ and the SNR. Vector  $\by$ was then calculated at 
$n$ observation points as $\by = \bq + \sig \,  \boeta$ where $\boeta \in \RR^n$ is a standard normal vector. 
Finally, we ran simulations for   two noise levels:  SNR = 3  and  SNR = 5. 
% respectively with different values of sample points as mentioned earlier.
% What do you mean by "different values of sample points"?

We compared the estimators of $\bof$ based on  random dictionaries with the estimator of $\bof$ 
based on the Singular Value Decomposition (SVD). 
For our simulations we have created  three  different   $n \times p$ random dictionaries   with $p =5000$:   
(a) two purely random dictionaries with, respectively, the  i.i.d. standard Gaussian entries and the  i.i.d. sparse Bernoulli entries;
(b) the fusion of the fixed dictionary and the random dictionary that follows case c) in Lemma~\ref{lem1} with $\bD$ being the Haar dictionary. 
The  sparse Bernoulli variable is defined as 
\bes
\bX=\left \{
\begin{array}{ll}
      -\sqrt{\frac{3}{n}} & \text{with probability} \ \frac{1}{6} \\
      0 & \text{with probability} \  \frac{2}{3} \\
     \sqrt{\frac{3}{n}} & \text{with probability} \ \frac{1}{6} \\
\end{array} 
\right. 
\ees
For creating the fusion dictionary, we first generated the orthogonal matrix of the Haar wavelet transform $\bD$ using {\tt MakeWavelet} function,
so that $m=n$ and $c=1$. Then we obtained the dictionary $\bPhi$ following part c) of the Lemma \ref{lem1} using the $n \times p$ matrix $\bW$ 
 with the i.i.d. normal entries.

% In order to generate operator $\bQ$ we sampled functions $f$ and $g$ on a fine grid.    
We obtained matrix $\bPsi$ of the inverse images as the numerical  solution of the exact equation $\bQ^T \bPsi = \bPhi$ and  
calculated  vector $\bz$ with elements \fr{bz}.
For the sake of obtaining a solution of optimization problem \fr{las_sol}, we used function {\tt LassoWeighted} 
in  SPAMS MatLab toolbox (see \cite{spams}). % Mairal  (2014)).

In order to   evaluate the value of the lasso parameter $\af$, we calculated $\af_{\max}$ 
as the value of $\af$ that guarantees that all coefficients  in the model vanish.
We created a grid of the values of $\af_k = \af_{\max}*k/N$, $k=1, \cdots, N$, with $N = 200$. 
As a result, we obtained a collection of estimators $\hbte = \hbte(\af_k)$. 
For the purpose of choosing the most appropriate  value of $k$,  
we estimated $\af$ as $\hat{\af} = \af_{\hk}$ in two ways: one using the oracle value of $\af$ 
and another using the estimated value of $\af$. We found oracle  value of $\alpha$ as
$\af_{oracle} = \af_{\max}*\hk_{oracle}/N$ using the value   $\hk_{oracle}$ 
that guarantees the most accurate estimator of $\bof$:
\bes % \label{ideal}
\hk_{oracle} = \arg\min_k   \| \bof- \bPhi \hbte(\af_k) \|_2 .
\ees
Since the vector $\bof$ is unavailable in real life, we find the  estimated value 
$\hat{\af}_{est}= \af_{\max}*\hk_{est}/N$ of $\af$ using 
\bes % \label{cv}
\hk_{est} = \arg\min_k \lfi  \frac{1}{n} \| \by - \hbq(\af_k) \|^2_2 + 2 \sig^2 n^{-1} \hat{p}_k \rfi,
\ees  
where $\hbq(\af_k) = \bQ \bPhi \widehat{\btheta} (\af_k)$ is the estimator of $\bq$ based 
on the lasso estimator obtained with the parameter $\af_k$ and  $\hat{p}_k$ is the number of nonzero components of $\hbte(\af_k)$. 

%%%%%%%%%%%%%%%%%%%%%%%%%%%%%%%%%%%%%%%%%%%%%%%%%%%%%%%%%%%%%%%%%%%%%%%=

 \begin{figure} %[ht]
 \[\includegraphics[height=4.3cm]{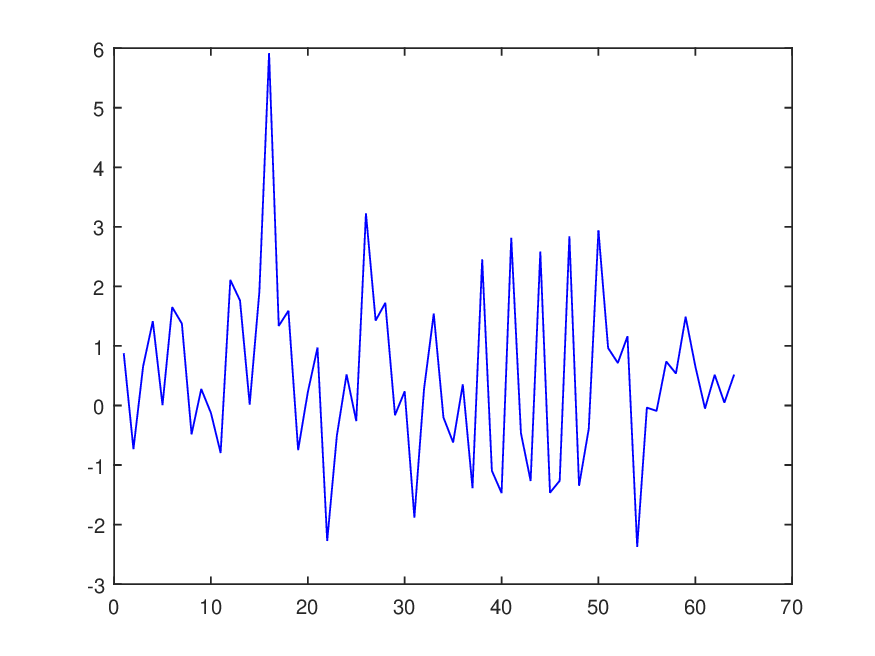} %\hspace{2mm}  
 \includegraphics[height=4.3cm]{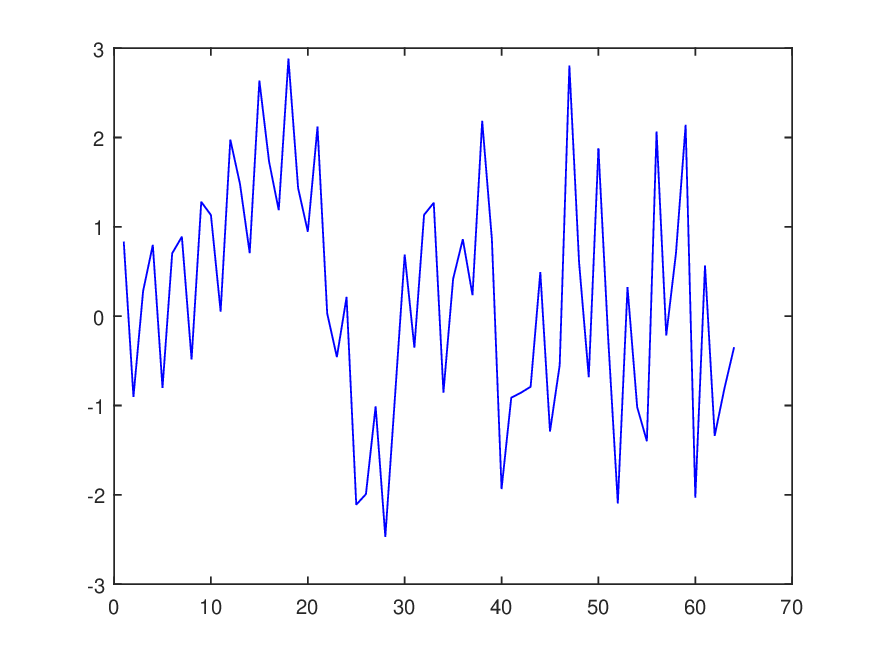} %\hspace{2mm}
 \includegraphics[height=4.3cm]{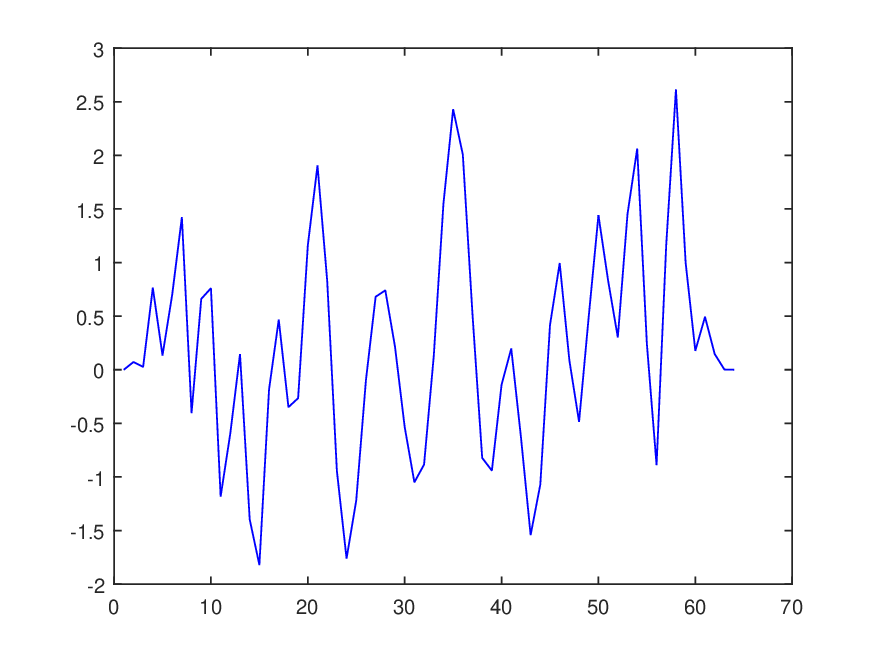} \]
 \caption{ Test signals {\tt WernerSorrows} (left), {\tt MishMash} (middle)
and {\tt Chirps} (right) with $n=64$.
\label{fig1}}
\end{figure}

%%%%%%%%%%%%%%%%%%%%%%%%%%%%%%%%%%%%%%%%%%%%%%%%%%%%%%%%%%%%%%%%%%%%%%%

We compared the estimators $\hbf_{RN}$, $\hbf_{RB}$, $\hbf_{RH}$  of $\bof$ based, respectively, on 
Gaussian, Bernoulli and Haar fusion random dictionaries described above  with  $\hbf_{SVD}$, 
the estimator based on the singular value decomposition (SVD).  
Initially we considered wavelet estimator of $\bof$ using Daubechies wavelet of order 8, 
but we discarded it due to its poor performance with respect to the estimators considered for comparison.
For finding $\hbf_{SVD}$, we used the oracle number $K_{oracle}$  of eigenbasis functions. 
We obtained  $K_{oracle}$  as the number of eigenbasis functions  that minimizes the difference between  $\hbf_{SVD}$ and the true function
$\bof$ which is unavailable in a real life setting.

Table \ref{table1} below compares the accuracies  of the estimators based on random dictionaries with the SVD estimator. 
Precision of an estimator $\hbf$  is measured by  $R(\hbf) =    n^{-1/2}\,  \|\hbf - \bof\|_2$,
the estimated $L^2$-norm of the difference between the estimator $\hbf$ and the true vector $\bof$
averaged over 50 simulation runs  (with the  standard deviations listed  in parentheses). 
For all the three estimator based on random dictionaries, we report the errors with both the oracle and the estimated values of $\af$,
$\hbf_{RN}^{oracle}$, $\hbf_{RB}^{oracle}$, $\hbf_{RH}^{oracle}$ and $\hbf_{RN}^{est}$, $\hbf_{RB}^{est}$, $\hbf_{RH}^{est}$,
respectively. 
We carried out simulations with three types of test functions {\tt WernerSorrows},  {\tt MishMash} and {\tt Chirps}.
The test signals are presented in Figure \ref{fig1}.

%%%%%%%%%%%%%%%%%%%%%%%%%%%%%%%%%%%%%%%%%%%%%%%%%%%%%%%%%%%%%%%%%%%%%%%%%%%%%%%%%%%%%%%%%%%%%%%% 

  \renewcommand{\baselinestretch}{1}

  \begin{table} %[ht]
\begin{center}
\begin{tabular}{|l| c |c | c | c | c | c | c | }
  % \multicolumn{4}{ c }{{\sc  averaged over 50 simulation runs  }}     \\
 \multicolumn{8}{ c }{ $WernerSorrows$  }\\
\hline
 &     $\hbf^{oracle}_{RN}$   &   $\hbf^{est}_{RN}$   &  $\hbf^{oracle}_{RB}$   &   $\hbf^{est}_{RB}$   &    
$\hbf^{oracle}_{RH}$   &   $\hbf^{est}_{RH}$   &   $\hbf_{SVD}$   \\
\hline
\hline
$n = 32$,   & 0.4910      & 0.5127   & 0.4910  & 0.5125 & 0.4837  & 0.4981 & 0.5155  \\
 $SNR=3$  & (0.0715)  & (0.0767) & (0.0736) & (0.0803) & (0.0724) & (0.0716) & (0.0721)\\
\hline
 $n = 32$,  & 0.3838       & 0.3956     & 0.3810    &  0.3966   & 0.3801    & 0.3917  & 0.4020   \\
$SNR=5$    & (0.0642)  & (0.0632) & (0.0654) & (0.0637) & (0.0622) & (0.0631) & (0.0632)\\
 \hline
$n = 64$,      & 0.5818        & 0.6074     & 0.5874    & 0.6112  & 0.5759  & 0.5929   & 0.6114  \\
$SNR=3$   & (0.0511)  & ( 0.0541) & (0.0474) & (0.0537) & (0.0461) & (0.0510) & (0.0577)\\
\hline
 $n = 64$,     & 0.3152       & 0.3208    & 0.3166   & 0.3210   &  0.3160   & 0.3199  & 0.3254  \\
 $SNR=5$  & (0.0351)  & (0.0350) & (0.0360) & (0.0364) & (0.0349) & ( 0.0359) & (0.0359)\\
 \hline
 \hline
$n = 128$,     & 0.5761      & 0.6072   & 0.5780    & 0.6082  & 0.5744  & 0.5901  & 0.6218  \\
$SNR=3$   & (0.0406)  & (0.0431 ) & (0.0394) & (0.0442) & (0.0393) & (0.0417) & (0.0440)\\
\hline
$n = 128$,     & 0.3730   & 0.3807   & 0.3717  & 0.3801  & 0.3717   &  0.3757   & 0.3850   \\
$SNR=5$   & (0.0251)  & (0.0267) & (0.0244) & (0.0268) & (0.0250) & (0.0259) & (0.0271)\\ 
\hline 
\hline
  \multicolumn{8}{ c }{ }\\
  \multicolumn{8}{ c }{$MishMash$}\\
 \hline
  $n = 32$,   & 0.5409       & 0.5628   & 0.5432   & 0.5631  & 0.5406   & 0.5562   & 0.5977  \\
 $SNR=3$  & (0.0811)  & (0.0871) & (0.0800) & (0.0854) & (0.0848) & (0.0846) & (0.0874)\\
 \hline
 $n = 32$,  & 0.3523   &  0.3595   &0.3537   &  0.3613  & 0.3532   & 0.3592  & 0.3657  \\
 $SNR=5$    & (0.0588)  & (0.0578) & (0.0561) & (0.0566) & ( 0.0554) & (0.0560) & (0.0549)\\
 \hline
 $n = 64$,      & 0.6131   &  0.6400   & 0.6158   & 0.6391  & 0.6145  & 0.6352   & 0.6599  \\
 $SNR=3$   & (0.0520)  & (0.0617) & (0.0552) & (0.0595) & (0.0541) & (0.0606) & (0.0621)\\
 \hline
$n = 64$,     &  0.3039  &  0.3086  & 0.3039  & 0.3083  &  0.3024  & 0.3067  & 0.3107  \\
 $SNR=5$  & (0.0298)  & (0.0315) & (0.0283) & (0.0300) & (0.0293) & (0.0294) & (0.0297)\\
 \hline
 \hline
 $n = 128$,     & 0.5112      & 0.5252   & 0.5106    & 0.5246  & 0.5110  & 0.5185  & 0.5410 \\
 $SNR=3$   & (0.0393)  & (0.0409) & (0.0381) & (0.0412) & (0.0394) & (0.0402) & (0.0415)\\
 \hline
 $n = 128$,     & 0.3385   & 0.3430   & 0.3380  & 0.3429  & 0.3383   &  0.3410   & 0.3460  \\
 $SNR=5$   & (0.0244)  & (0.0240) & (0.0241) & (0.0242) & (0.0237) & ( 0.0235) & (0.0245)\\ 
 \hline \hline
  \multicolumn{8}{ c }{ }\\
 \multicolumn{8}{ c }{$Chirps$}\\
 \hline
 $n = 32$,   & 0.4300   & 0.4430   & 0.4320  & 0.4431  & 0.4307  & 0.4443 & 0.4607 \\
 $SNR=3$  & (0.0630)  & (0.0646) & (0.0602) & (0.0598) & (0.0613) & (0.0623) & (0.0607)\\
 \hline
 $n = 32$,  & 0.2872  &  0.2977   & 0.2880  &  0.2960  & 0.2871 & 0.2951 & 0.3018\\
 $SNR=5$    & (0.0461)  & (0.0459) & (0.0448) & ( 0.0448) & (0.0457) & (0.0444) & (0.0443)\\
 \hline
 $n = 64$,      & 0.3979   &  0.4117  & 0.3983  & 0.4121 & 0.3986  & 0.4066  & 0.4297 \\
 $SNR=3$   & (0.0369)  & (0.0404) & (0.0393) & (0.0414) & (0.0406) & (0.0404) & (0.0391)\\
 \hline
 $n = 64$,     &  0.2735   & 0.2782    & 0.2733  &  0.2778  &   0.2723 &  0.2767  & 0.2789 \\
 $SNR=5$  & (0.0314)  & (0.0345) & (0.0327) & (0.0351) & (0.0320) & (0.0341) & ( 0.0347)\\
 \hline
 \hline
 $n = 128$,     & 0.3878   & 0.3955   & 0.3867   & 0.3953 & 0.3868 & 0.3953 & 0.4069  \\
 $SNR=3$   & (0.0288)  & ( 0.0310) & (0.0282) & (0.0289) & ( 0.0289) & ( 0.0285) & (0.0293)\\
 \hline
 $n = 128$,     & 0.2428  & 0.2460   & 0.2431 & 0.2456  & 0.2428 &  0.2460  &  0.2473 \\
 $SNR=5$   & (0.0182)  & (0.0179) & (0.0180) & (0.0179) & (0.0180) & (0.0181) & (0.0173)\\ 
 \hline 
\hline
\end{tabular}
\end{center}
\caption{ The average values of the errors  $R(\hbf)$ evaluated over 50 simulation runs of the estimators for 
various test signals (with the standard deviations of the errors  listed in the parentheses).
} \label{table1}
\end{table}

From    Table \ref{table1}  it follows that all the random dictionary based estimators are 
more accurate than the  SVD estimator. The advantage of $\hbf_{RN}^{oracle}$, $\hbf_{RB}^{oracle}$ and  $\hbf_{RH}^{oracle}$ 
over $\hbf_{SVD}$ is more significant than that of 
 $\hbf_{RN}^{est}$, $\hbf_{RB}^{est}$ and $\hbf_{RH}^{est}$ since the latter estimators 
loose accuracy because of    suboptimal choices of the parameter $\alpha$. Nevertheless, 
in majority of cases, they still exhibit better precision than $\hbf_{SVD}$ although this is not entirely
fair comparison since $\hbf_{SVD}$ is based on the oracle choice  of parameter $K$.
This is due to the fact that large random dictionaries provide a more sparse representation of $\bof$.

%%%%%%%%%%%%%%%%%%%%%%%%%%%%%%%%%%%%%%%%%%%%%%%%%%%%%%%%%%%%%%%%%%%%%%%%%%%%%%%%%%%%%%%%%%%%%%%%%%%%%%%%%%%%%%%%%%%%
 
\section{Discussion}
\label{sec:discussion}

In the present paper we provided a new approach for the solution of a general ill-posed linear inverse problem. 
The underlying idea is to use lasso technique for estimating the function of interest 
by representing it as a sparse linear combination of elements of a random overcomplete dictionary. 
The advantage of choosing a random dictionary over any other overcomplete dictionary is that 
one can construct it in such a way that it satisfies restricted isometry condition with a high probability
and, therefore, ensures that the compatibility condition (which guarantees fast convergence rates for lasso) also holds.

We provide theoretical justification for application of the lasso technique with the random dictionaries 
for solution of the linear inverse problems. We also support our theory by the simulation studies
which show that the proposed estimators have higher accuracy than the SVD estimators in spite 
of the fact that the SVD estimators are based on the oracle parameter choices. For this reason,
the advantage of the random dictionary based estimators is more significant when they are likewise constructed 
with the oracle choices of parameter $\alpha$. In fact, this is the part where our method has some room for improvement:  
since our procedure for estimating parameter $\alpha$ is rather elementary, it can   be fine-tuned using more 
advanced techniques.

% random dictionaries satisfy compatibility conditions for dictionaries to recover function of interest efficiently. 
% Also, random dictionaries satisfies restricted isometry condition with a high probability. 
% We provided the theoretical proof of better precision and supported it by simulation studies over 
% signals like WernerSorrows, MishMash, Blocks and Chirps using different sample points and 
% different amount of noise. The results in the simulation section demonstrates the reliability 
% of our methods over one of the other most commonly used methods, SVD. 

%%%%%%%%%%%%%%%%%%%%%%%%%%%%%%%%%%%%%%%%%%%%%%%%%%%%%%%%%%%%%%%%%%%%%%%%%%%%%%%%%%%%%%%%%%%%%%%%%%%%%%%%%%%%%%%%%%%%

 \section*{Acknowledgements}

Marianna Pensky  and Pawan Gupta were  partially supported by National Science Foundation
(NSF), grants   DMS-1106564 and DMS-1407475.

%%%%%%%%%%%%%%%%%%%%%%%%%%%%%%%%%%%%%%%%%%%%%%%%%%%%%%%%%%%%%%%%%%%%%%%%%%%%%%%%%%%%%%%%%%%%%%%%%%%%%%%%%%%%%%%%%%%%

\section{Proofs} 
\label{sec:proofs}
\setcounter{equation}{0}

{\bf Proof of Theorem \ref{th:randLasso}.  }
The beginning of the proof is similar to the proof of Lemma 2 in \cite{pensky}.
However, for completeness, we provide the complete proof here.

Let $\bte$ be the true parameter vector, so that $\bof = \bof_{\bte} = \bPhi \bte$. 
Denote $\bzeta = \bPsi^T \boeta$. Then, it is easy to check that 
\bes
\bPhi^T (\bz - \bof) = \bPsi^T(\by - \bQ \bof) = \sig \bzeta.
\ees
Following \cite{arnak}, by K-K-T condition, we derive that for any $\bt \in \RR^p$
\beqns
\hbte^T \bPhi^T\, (\bz - \bPhi \hbte) & = & \af \sumjp \nuj |\hte_j| \\
\bt^T  \bPhi^T\, (\bz - \bPhi \hbte) & \leq & \af \sumjp \nuj |t_j|, 
\eeqns
so that, subtracting the first line from the second,  we obtain
\be \label{main_ineq}
(\bPhi \hbte - \bPhi \bt)^T (\bPhi \hbte - \bz) \leq \af \sumjp \nuj (|t_j| - |\hte_j|).
\ee 
Then, \fr{main_ineq} yields
%\bes 
$
(\bPhi \hbte - \bPhi \bt)^T   (\bPhi \hbte - \bPhi \bte) \leq \sig  (\hbte - \bt)^T \bzeta + \af \sumjp \nuj (|t_j| - |\hte_j|).
$
%\ees
Since for any $\bu, \bv \in \RR^p$ one has 
%\beqns
$
\bv^T   \bu = \frac{1}{2} \lkv \|\bv \|^2 + \|\bu\|^2 - \|\bv - \bu\|^2 \rkv,
$
%\eeqns
choosing $\bv = \bPhi \hbte - \bPhi \bt$ and $\bu = \bPhi \hbte - \bPhi \bte$  
for any $\bt  \in \RR^p$ obtain
\be \label{ineq1} 
\| \bof_{\hbte} - \bof \|^2 + \|\bPhi   (\hbte -   \bt)\|^2 \leq \| \bof_{\bt} - \bof \|^2  + 
2 \sigma   (\hbte - \bt)^T \bzeta + 2 \af \sumjp \nuj (|t_j| - |\hte_j|).
\ee 
By definition of  $\bzeta$, for any $j = 1, \cdots, p$, one has $\zeta_j \sim \calN (0, \nuj^2)$.
Hence, on the set 
\be \label{eq:setOm0} 
\Om_0 = \lfi \om: \max_{1 \leq j \leq p} (\nuj^{-1} |\zeta_j|) \leq   \sqrt{2 (\tau +1) \log p} \rfi \quad \mbox{with}
\quad \PP(\Om_0) \geq 1 -  2 p^{- \tau}
\ee 
one obtains 
%\bes 
$
|(\hbte - \bt)^T \bzeta| \leq \sqrt{2 (\tau +1) \log p}\, \sumjp \nuj \, |\hte_j - t_j| = \af_0 \, \sumjp \nuj \, |\hte_j - t_j|. 
$
%\ees
Combining the last inequality with \fr{ineq1} 
obtain that, for any $\af >0$, on the set $\Om_0$,
\be \label{ineq2} 
\| \bof_{\hbte} - \bof \|^2 + \|\bPhi   (\hbte -   \bt)\|^2 \leq \| \bof_{\bt} - \bof \|^2    + 
2  \af \sumjp \nuj (|t_j| - |\hte_j|) + 2 \alfo \sumjp \nuj \, |\hte_j - t_j|.
\ee 
Denote   
$\bOm_1 = \lfi \om:\ \lamin(2s; \bPhi)   \geq    1-\delta \rfi$, 
$\bOm_2 = \lfi \om:\ | J_{*}| \leq s \rfi$  and 
$\bOm_3 = \lfi \om:\  \|\hbte\|_0 \leq s \rfi$.
Choose $\bt$ such that $\bof_{\bt} = \proj_{\calL_{J_{*}}} \bof = \bof_{\calL_{J_{*}}}$ and 
note that $t_j = 0$ for $j \in J_{*}^c$.
Then, due to $\af \geq \af_0$ and $||\hte_j - t_j| \leq |\hte_j| + |t_j|$,
obtain 
\be \label{ineq3} 
\| \bof_{\hbte} - \bof \|^2 + \|\bPhi   (\hbte -   \bt)\|^2 \leq \| \bof_{\bt} - \bof_{\calL_{J_{*}}} \|^2    + 
4 \af \ \sumJ \nuj \, |\hte_j - t_j|.
\ee 

Consider the set $\bOm = \bOm_0 \cap \bOm_1\cap \bOm_2\cap \bOm_3$ and note that 
$\PP (\bOm) \geq 1 - 2 p^{-\tau} - \eps_1 - \eps_2 - \eps_3.$
If $\om \in \bOm$, then $\|\hbte - \bt\|_0 \leq 2s$ and, hence,
\bes
4 \af \ \sumJ \nuj \, |\hte_j - t_j| \leq 
4 \af \lkr \sumJ \nuj^2 \rkr^{1/2} \, \frac{\|\bPhi   (\hbte -   \bt)\|}{\lamin (2s; \bPhi)}
\leq  \|\bPhi   (\hbte -   \bt)\|^2 + \frac{4 \af^2}{(1 - \del)^2} \, \sumJ \nuj^2.
\ees
Plugging the last inequality into \fr{ineq3} and recalling the definition of $J_{*}$, we derive 
\fr{newfast1}.
\\

%%%%%%%%%%%%%%%%%%%%%%%%%%%%%%%%%%%%%%%%%%%%%%%%%%%%%%%%%%%%%%%%%%%%%%%%%%%%%%%%%%%%%%%%%%%%%%%%%

\medskip

{\bf Proof of Lemma \ref{lem1}.  }
In cases a) and b), $\lamin (m; \bPhi)\geq 1-\del$ 
% validity of \fr{lammin} 
is ensured by Theorem~5.65 of {Vershynin  \cite{vershynin}.
In case c),  note that entries of matrix $\bPhi$ are uncorrelated and, hence, are independent Gaussian variables due to 
 \bes
 \Cov\left(\bPhi_{ik}\cdot \bPhi_{jl}\right)  =  \frac{1}{c^2}\sum_{r_1=1}^{m} \sum_{r_2=1}^{m}\bD_{ir_1}\bD_{jr_2}I(r_1=r_2)I(k=l)
= I(i=j)I(k=l). 
 \ees
Moreover, matrix $\bPhi$ has isotropic rows since  
 \bes
 \Cov\left(\bPhi_{ih}\cdot \bPhi_{jl}\right)  =  \frac{1}{c^2}\sum_{r_1=1}^{m} \sum_{r_2=1}^{m}\bD_{ir_1}\bD_{jr_2}I(r_1=r_2)I(h=l)
= I(i=j)I(h=l). 
 \ees
Therefore, 
$\lamin (m; \bPhi)\geq 1-\del$ by Theorem 5.65 of \cite{vershynin}.
\\

%%%%%%%%%%%%%%%%%%%%%%%%%%%%%%%%%%%%%%%%%%%%%%%%%%%%%%%%%%%%%%%%%%%%%%%%%%%%%%%%%%%%%%%%%%%%%%%%%%%%%%%%%%%%%%%%%%%%55

\end{document}